\documentclass{amsart}
\usepackage{amsmath}
\usepackage{amsthm}
\usepackage{faktor}
\usepackage{amssymb}
\usepackage{tikz-cd}
\usepackage{courier}
\usepackage{amsfonts}
\usepackage{multicol}
\usepackage[mathscr]{euscript}
\usepackage{dirtytalk}
\usepackage{graphicx}
\graphicspath{{./images/}}
\newcommand{\interior}[1]{%
  {\kern0pt#1}^{\mathrm{o}}%
  }
\newtheorem{theorem}{Theorem}[section]
\newtheorem{lemma}[theorem]{Lemma}

\newtheorem{definition}{Definition}[section]
\newtheorem{conjecture}[theorem]{Conjecture}

\newtheorem*{theorem*}{Theorem}

\begin{document}
\title{\texttt{Plat closures of spherical braids in $\mathbb{R}P^3$}}
\author{\texttt{\textbf{Rama Mishra, Visakh Narayanan\\
Indian Institute of Science Education and Research, Pune.}}}
\maketitle
\begin{abstract}
We define plat closure for spherical braids  to obtain links in $\mathbb{R}P^3$ and prove that all links in $\mathbb{R}P^3$ can be realized in this manner. Given a spherical braid $\beta$ of $2n$ strands  in $\mathbb{R}P^3$ we associate a permutation $h_{\beta}$ on $n$ elements called \textit{residual permutation}. We prove that  the number of components of the plat closure link of a spherical braid $\beta$ is same as the number of disjoint cycles in $h_{\beta}$. We also present a set of moves on spherical braids in the same spirit as the classical Markov moves on braids. The completeness of this set of moves to capture the entire isotopy classes of the plat closure links is still to be explored.  
\end{abstract}
\section{\texttt{\textbf{Introduction}}}
The discovery of quantum invariants has revolutionized the study of low dimensional topology. For example, in case of classical knots, the long standing problem of the trefoil and its mirror image being not equivalent, was solved by the Jones polynomial, one of the first quantum invariants. The discovery and computation of these invariants have been facilitated by the intimate relationship between Artin's braid group and classical knots. \\

As knot theory grows into a basic subject in low dimensional topology, it becomes important to understand the different \say{ramifications} of it such as studying knots inside other three-manifolds. The real projective 3-space, being one of the closest cousins of $S^3$, it is natural to consider the knot theory it admits as the next candidate. But the important lesson we learn from classical knot theory is that, a knot and the space around it cannot be separated from each other. Their topological features are closely connected. Hence the complexity of the ambient manifold raises natural challenges in its knot theory.\\

In this paper we introduce a braid theory for the knots and links in real projective 3-space. Joan Birman \cite{birman} introduced the notion of braid group of an arbitrary manifold. In that setup, the standard Artin's braid group is the braid group of the plane. The braid group of the 2-sphere is also of particular importance. We will refer to the elements of this group as \say{spherical braids} and the elements of Artin's braid group as \say{classical braids}. The $n$-string braid group of any manifold may be described as the group of motions of $n$ special points in it. Birman then goes on to define the concept of plats in $S^3$ \cite{birman} as a different \say{closure} of braids than the standard closure \cite{alex}, and proves that every classical link is isotopic to a plat. In this paper, we develop the notion of plats in $\mathbb{R}P^3$ using spherical braids and show that \textit{every link in $\mathbb{R}P^3$ is isotopic to a plat (Theorem 3.2)} defined in this way. We call it as \say{plat closure} of a spherical braid. For the sake of completeness, clear definitions of these are provided in the Section 2.   \\

Another beautiful feature of a classical braid is the permutation that it defines. The number of disjoint cycles in the permutation of the classical braid is equal to the number of components in the closure link. In Section 5, we define a permutaion which we refer to as \say{residual permutation}, assosiated to a spherical braid in $\mathbb{R}P^3$. We prove that the number of disjoint cycles in the residual permutation of a spherical braid matches with the number of components in its plat closure (Theorem 5.1). Towards the end of this section we discuss a theorem which characterizes when will the plat closure of a braid be affine.  \\

J.W. Alexander \cite{alex}, proved that every classical link is isotopic to the closure of some classical braid. This representation is not unique as many braids may close to give isotopic links. Andrei Markov defined a set of moves on classical braids, which are now known as \say{Markov moves}, and showed that two braids can have isotopic closures if and only if they are related by finitely many Markov moves \cite{birman}. In Section 6, we present some moves on spherical braids, in the spirit of classical Markov moves. We call them $M$-moves and two braids which are related by a finite sequence of these moves are said to be $M$-equivalent. We conjecture that two plat closure links can be isotopic if and only if their braids are $M$-equivalent. \\
 
\textbf{Organization of the paper:} Section 2, provides definitions and terminologies which will be used throughout the paper. In Section 3, we include discussions on plat closures of spherical braids in $\mathbb{R}P^3$ and prove that every link in $\mathbb{R}P^3$ is isotopic to plat closure of some sherical braid. In Section 4, we study some properties of the braid group of the 2-sphere and give a convenient presentation for the group $B_n(S^2)$. In Section 5 we define the notion of a \say{residual permutation} of a spherical braid and prove that the number of components in the closure link is same as the number of disjoint cycles in this permutation. Section 6 introduces a set of moves on spherical braids (M-moves) which possibly can capture the isotopy class of the plat closure link of a given spherical braid.

\section{Spherical braids}
Let $\Sigma$ be a manifold of arbitrary dimension and let $n$ be a positive integer. Consider a set of $n$ distinct points $X:=\{p_1,p_2,...,p_n\}\subset \Sigma$. These are thought of as special $n$ points in space. Consider an isotopy,\\
\begin{equation*}
F:\Sigma\times [0,1]\to \Sigma,
\end{equation*}\\
such that, $F_0=Id_\Sigma$ and $F_1$ is a diffeomorphism of $\Sigma$ mapping $X$ to itself. Then we may represent the isotopy in the space $\Sigma\times I$ by considering the map,\\
\begin{align*}
    \overline{F}:\Sigma\times I\to \Sigma\times I,\\
    (q,t)\mapsto (F_t(q),t)
\end{align*}\\
defined by it. Notice that the image of the set $X\times I$ under $\overline{F}$ is a collection of paths each one starting at some $(p_i,0)$ and ending at some $(p_j,1)$. \\
\begin{definition}
The topological pair $(\Sigma\times I,\overline{F}(X\times I))$ is refered to, as an \textbf{$n$ braid of $\Sigma$}.
\end{definition}

 Consider the case when $\Sigma=S^2$. Each $n$-braid of $S^2$ is represented by a set of $n$ paths in $S^2\times I$ such that each path starts at a point of $S^2\times \{0\}$ and ends at $S^2\times \{1\}$ and intersects each of the sections $S^2\times \{i\}$ at a unique point transvesally. We may interpret the strip $S^2\times I$ as a 2+1 dimensional spacetime where the $I$ direction represents the flow of time. If we do not allow points of $S^2$ to move back in time their world lines intersect each \say{spacelike} sphere in the strip, i.e., spheres of the form $S^2\times \{t\}$ at a unique point, just like the strings of braids. Or in other words, the projection map $f:S^2\times I\to I$, is monotonic when restricted to each of the strings. Let $\alpha$ and $\beta$ be two such motions of points in a set $X\subset S^2$. Then clearly we can define a new motion by composing them, that is performing $\beta$ after $\alpha$, which we will call $\alpha\beta$. As a braid, it is defined as the braid obtained by gluing the $S^2\times\{1\}$ of the strip containing $\alpha$ to the $S^2\times\{0\}$ of the strip containing $\beta$ matching the indices properly and then rescaling the newly formed strip. This defines a multiplication of braids. The set of isotopy classes of braids forms a group under this operation. We describe this group in detail in the Section 4.\\
 
\par Let $B$ be the 3-ball. We know that, $\partial (S^2\times I)\approx S^2\amalg S^2$. Notice that $S^3$ can be obtained by gluing boundaries of $B\amalg B$ and $S^2\times I$. Classically the plats in $S^3$ were constructed \cite{birman} by considering a spherical braid in this strip and certain simple tangles (which we discuss in the next section) in both the balls. We wish to discuss a generalization of this construction. Let $M$ denote the mapping cylinder of the canonical two sheeted covering map $S^2\to \mathbb{R}P^2$. Notice that, $\partial M\approx S^2$. By gluing the boundaries of $M\amalg B$ and $S^2\times I$ we can obtain a copy of $\mathbb{R}P^3$. We introduce a different set of tangles in $M$. When we say \say{braids in $\mathbb{R}P^3$}, we mean the braids in $S^2\times I$ region in some splitting of $\mathbb{R}P^3$ of this type. Now by considering a braid in $S^2\times I$ and gluing its boundary with the special tangles in $B$ and $M$ we can form a collection of linked knotted curves. We would refer to these as the \say{projective plat closure}, or simply, plat closure of the spherical braid in $\mathbb{R}P^3$.\\

Without loss of generality, we may assume that the special points are lying on an equatorial circle, $C$ of $S^2$. Then we can project every braid into an annulus, in a projective plane in $\mathbb{R}P^3$. Here we are using the same projection which were used in \cite{julia} and \cite{ramavis}. Thus we can represent such braids by a diagram drawn on the annulus. See to Figure 1. The diagram of a composition $\alpha\beta$ will appear as keeping the diagram of $\beta$ \say{inside} the diagram of $\alpha$.\\

\begin{figure}
\includegraphics{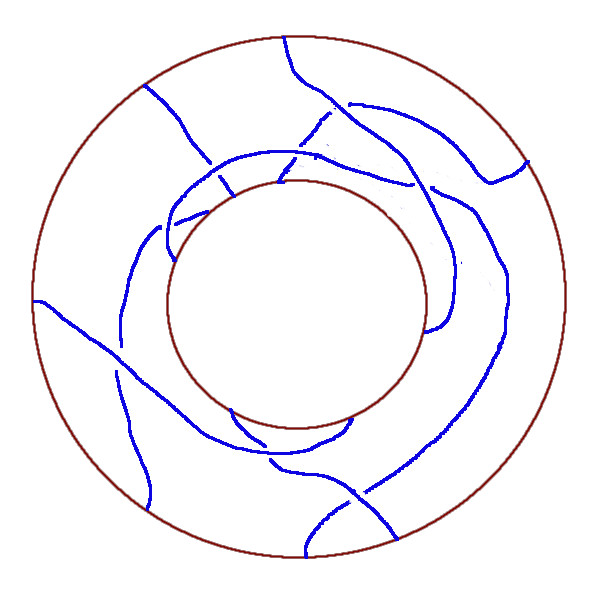}
\caption{Example of a braid on $S^2$}
\end{figure}
\par

\section{Plat closures in $\mathbb{R}P^3$.}
Here, we would consider \say{closures} of these braids in $\mathbb{R}P^3$. Choose an equator $C$ for a $\partial B$ and let $D$ be the flat disk in $B$ with boundary $C$. Let $A^n$ represent the tangle in $B$ formed by $n$ unknotted unlinked arcs neatly embedded in $B$ lying on $D$. We will call them as internal tangles. Refer to Figure 2. 

\begin{figure}
\includegraphics[scale=0.9]{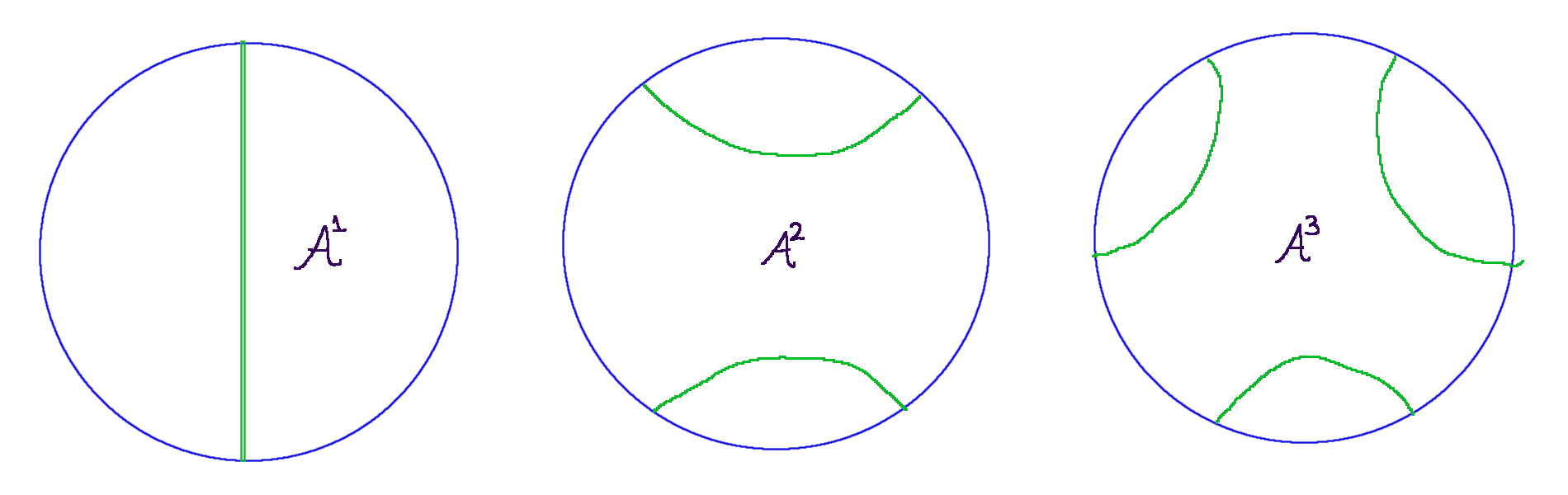}
\caption{\textbf{\texttt{Internal tangles}}}
\end{figure}
\par

Notice that the bourdary of the tangle $A^n$ is formed by $2n$ points on $\partial B$. We call certain special tangles in $M$, which appear in all links in $\mathbb{R}P^3$ as \say{residual tangles}. We may define them as follows: notice that $ H_1(M,\partial M)\approx \frac{\mathbb{Z}}{2\mathbb{Z}}$, every arc in $M$ with its two boundary points in $\partial M$, will represent a class in $H_1(M,\partial M)$. Residual tangles can be described as tangles formed by a collection of unknotted unlinked arcs each representing the $\bar{1}$ in $H_1(M,\partial M)$. We also require that all the arcs in the tangle are lying in a single flat Mobius band in $M$. Figure 3 demonstrates the first few examples. More properties of these are studied in \cite{ramavis}. 
\begin{figure}
\includegraphics{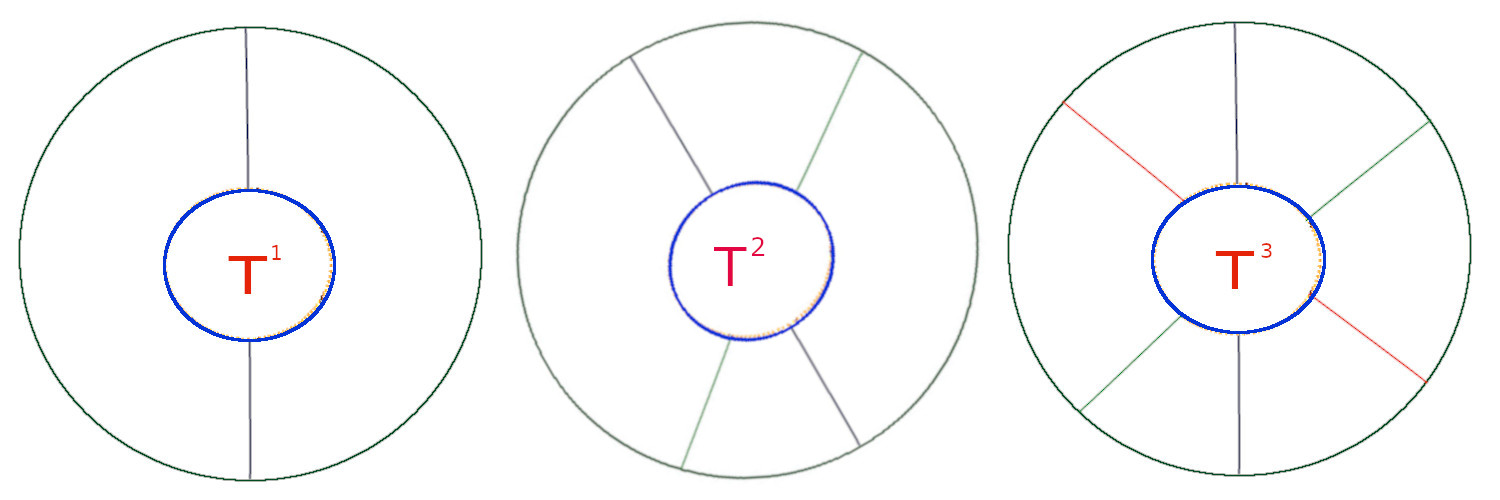}
\caption{\texttt{\textbf{Residual tangles}}}
\end{figure}

The boundary of the residual $n-$tangle is composed of $2n$ points. Now consider gluing the boundaries of $M$ with a residual $n$-tangle and $B$ with an internal $n-$tangle, using a diffeomorphism, $f:\partial M\to \partial B$, which sends $2n$ points on the boundary of $T^n$ to the $2n$ points on the boundary of $A^n$. By identifying $\partial M$ and $\partial B$ with $S^2$, we can find an isotopy $H:S^2\times I\to S^2$, of $f$ to the identity map of $S^2$. By representing the image of $\partial T^n$ under each of the maps, $h_t(x):=H(x,t)$ on the sphere $S^2\times \{t\}$ in $S^2 \times I$, we can obtain a braid in $S^2\times I$. Then the arcs in the internal tangle and residual tangle will join the boundary points of the braid. Thus we get a link in $\mathbb{R}P^3$. We will refer to this \say{closure} of braids as \textbf{projective plat closure}.

\begin{figure}
\includegraphics{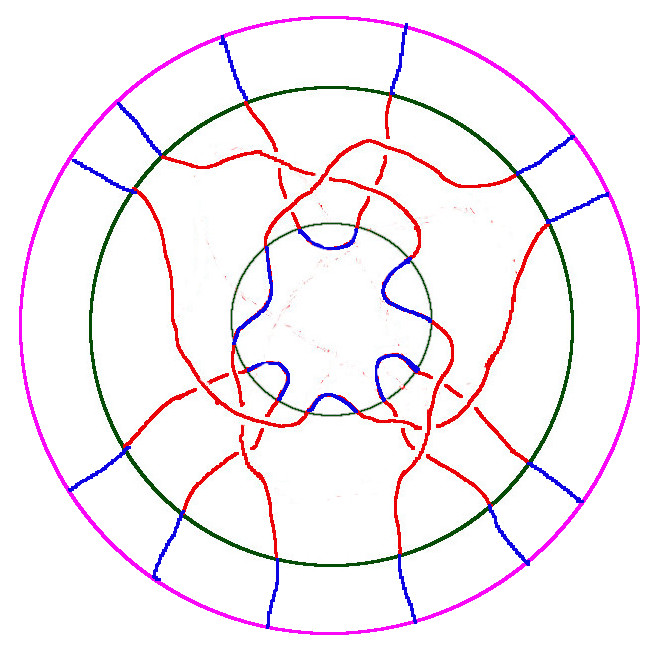}
\caption{A plat representation of affine trefoil}
\end{figure}

The following lemma will be used to prove our first theorem below. For the proof of the lemma, see Corollary 3.3 in \cite{ramavis}.

\begin{lemma}
Given any link $K$ in $\mathbb{R}P^3$, there exists a separating sphere, which will split $\mathbb{R}P^3$ into two pieces a ball $B$ and a mapping cylinder $M$ such that $K\cap M$ is a residual tangle. 
\end{lemma}

\begin{theorem}
Every link in $\mathbb{R}P^3$ is isotopic to the projective plat closure of a braid.
\end{theorem}

\textbf{Proof of the theorem:} Let $K$ be a link in $\mathbb{R}P^3$. Let  $S\subset \mathbb{R}P^3$ be the separating sphere provided by the lemma. Let $B$ and $M$ denote the ball and the mapping cylinder in the corresponding splitting respectively. The part of $K$ inside $M$ is already a residual tangle, say $T^n$. Let $B'\subset B$ be a smaller closed ball with the same center. See Figure 5. 
\begin{figure}
\includegraphics{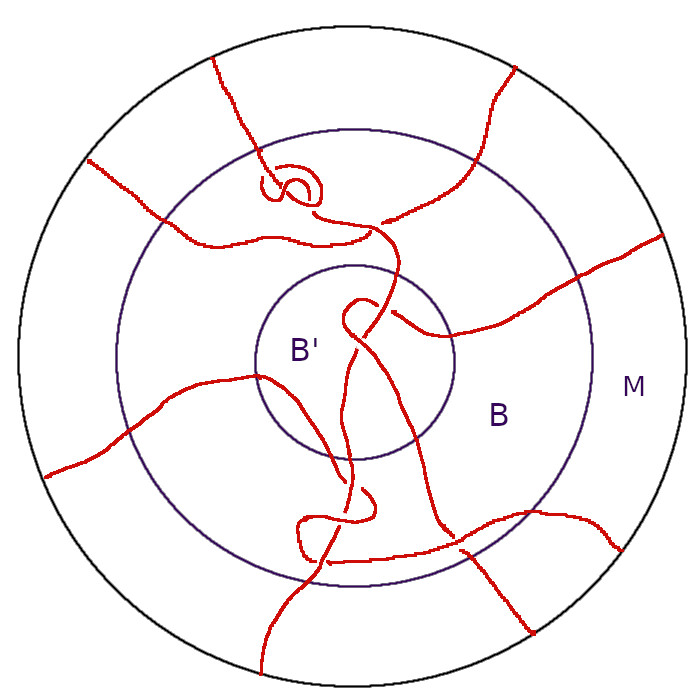}
\caption{A residual tangle in $M$ and a generic tangle inside $B$}
\end{figure}
Notice that the region outside $B'$ in $B$ homeomorphic to $S^2\times I$ with $S^2\times{0}$ mapped to $\partial B'$ and $S^2\times{1}$ mapped to $\partial B$. We will refer to this region as \say{\textbf{strip}} in what follows. \\
We can move the knot isotopically, so that all the crossings appear in the projection of the strip in the diagram.\\
\begin{figure}
\includegraphics{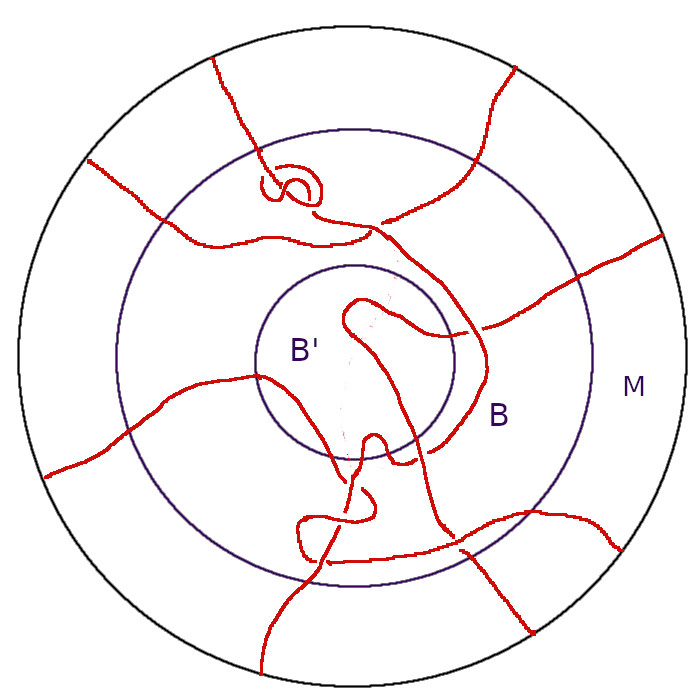}
\caption{After pushing all the crossings to the strip.}
\end{figure}
Now the tangle inside $B'$ is just a collection of untwisted unlinked arcs all whose boundaries are on $\partial B'$. The tangle inside the strip now has many arcs all of whose boundaries lie on $\partial B$ and $\partial B'$. The arcs may also be knotted. Refer to Figure 6. There are three types of arcs in the strip based on where their boundary points are placed. Both the boundary points of an arc may be on $\partial B'$. We will call them \say{type 1} arcs. The arcs with both the boundary points on $\partial B$ will be called \say{type 2} arcs. The arcs with one boundary point on $\partial B$ and another on $\partial B'$ will be called \say{type 3} arcs.   
\\
The projection,
\begin{equation*}
f: S^2\times I\to I
\end{equation*}
 on the strip can be restricted to the arcs. We will call this restriction, the \say{height  function} on the tangle. We shall denote this function also by $f$. Note that if $f$ has any point of inflection on an arc, the arc can be isotopically moved inorder to remove the inflection point and make $f$ monotonic locally. Hence, in what follows we will always assume that $f$ has no points of inflection on the tangle and extrema will mean either maxima or minima. Clearly on a type 1 arc, there exists atleast one maximum point for $f$. Similarly type 2 arcs has to have atleast one minimum point. If any of these arcs are knotted, they will contain more extrema points of $f$. Type 3 arcs which are not knotted may be isotopically moved so that $f$ is monotonic on the modified arc. Now all the extrema points may be removed from the strip by moving them inside $B'$. See Figure 7 for a demonstration. 
\begin{figure}
\includegraphics[scale=0.7]{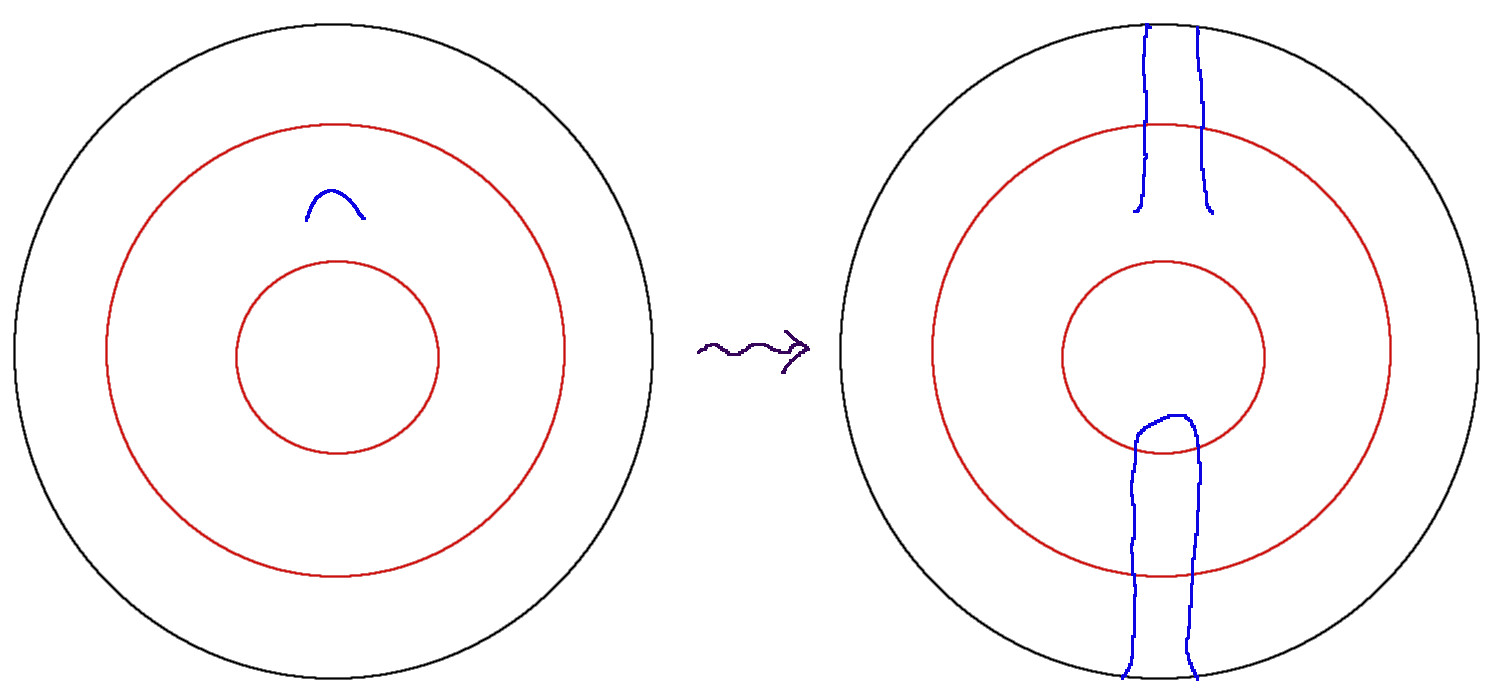}\\
\includegraphics[scale=0.7]{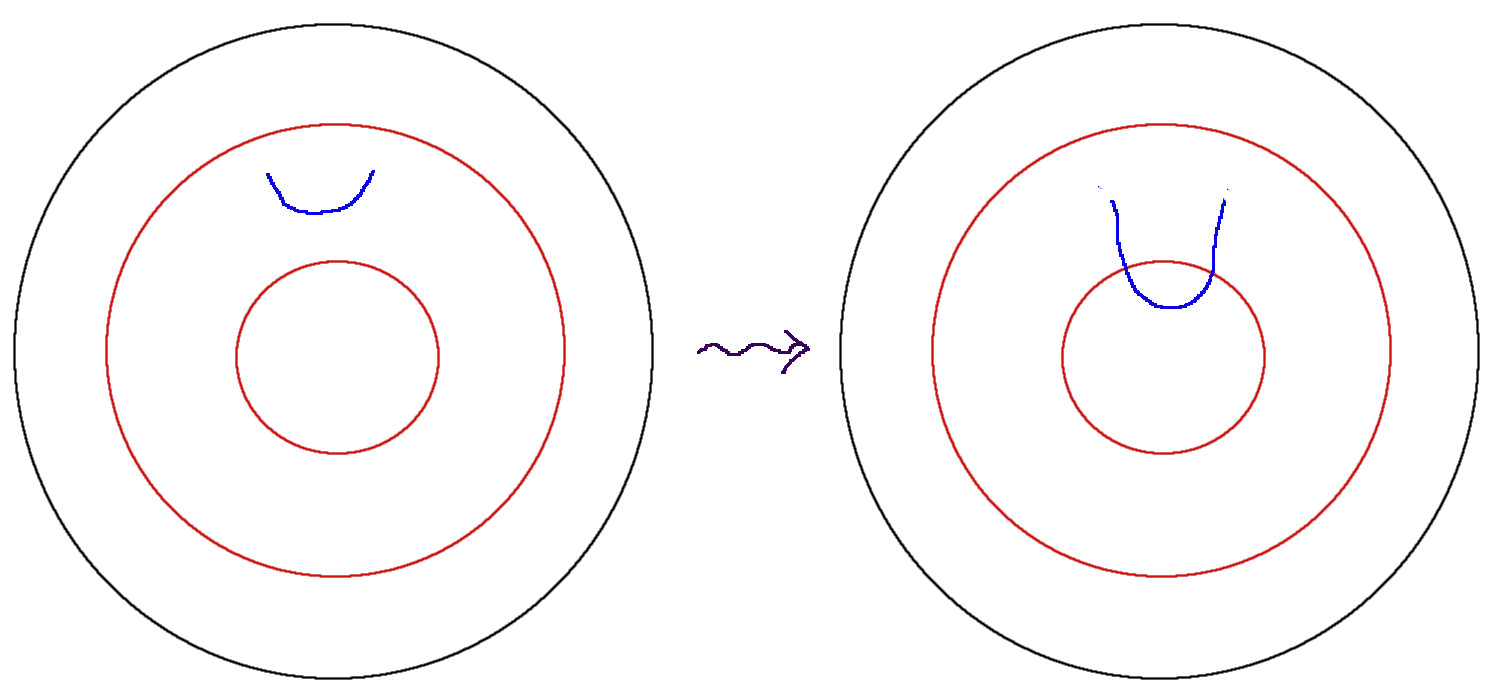}
\caption{Transfering the extremum points to $B'$.}
\end{figure}
It is easy to see that these operations can be done isotopically in $\mathbb{R}P^3$. And once all the extremum points are removed, $f$ will be monotonic on all the arcs in the strip. \\
Let $\gamma$ be an equator for $\partial B'$. It is easy to see that, we can isotopically move the link so that the tangle in $B'$ is an internal tangle with all the boundary points on $\gamma$. That is, all the boundary points on $\gamma$ are connected to their immediate neighbour. Refer to Figure 8.\\ 
\begin{figure}
\includegraphics{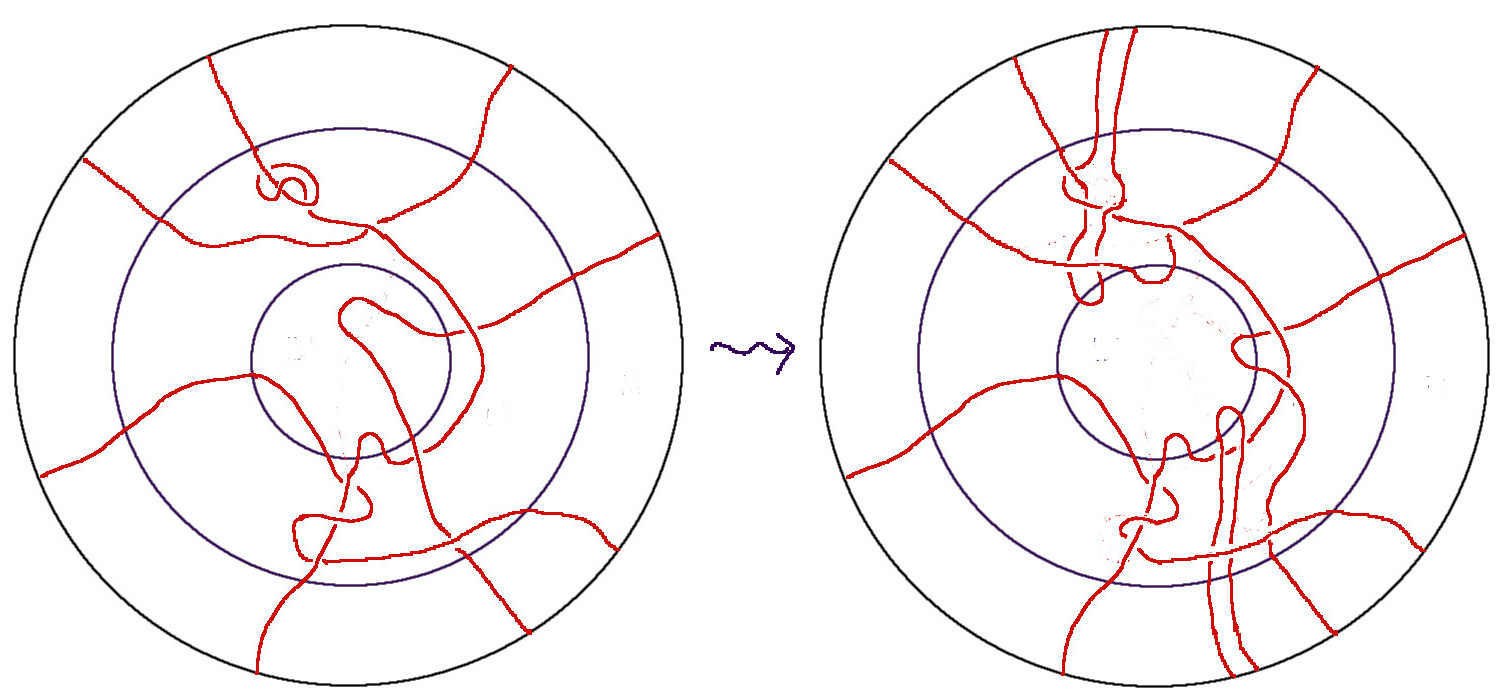}
\caption{An example.}
\end{figure}
\begin{figure}
\includegraphics[scale=1.2]{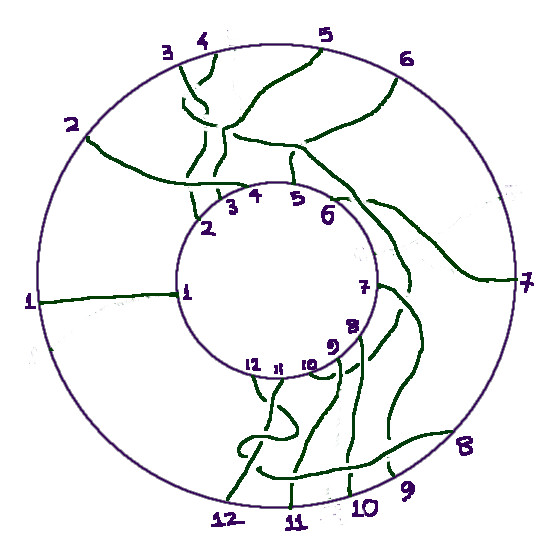}
\caption{The braid in the above example.}
\end{figure}
\\
Now it is easy to see that the tangle inside the strip is a braid. The residual tangle in $M$ and the internal tangle in $B'$ are \say{closing} this braid into a projective plat closure. Hence we are done. \qed

\begin{figure}
\includegraphics{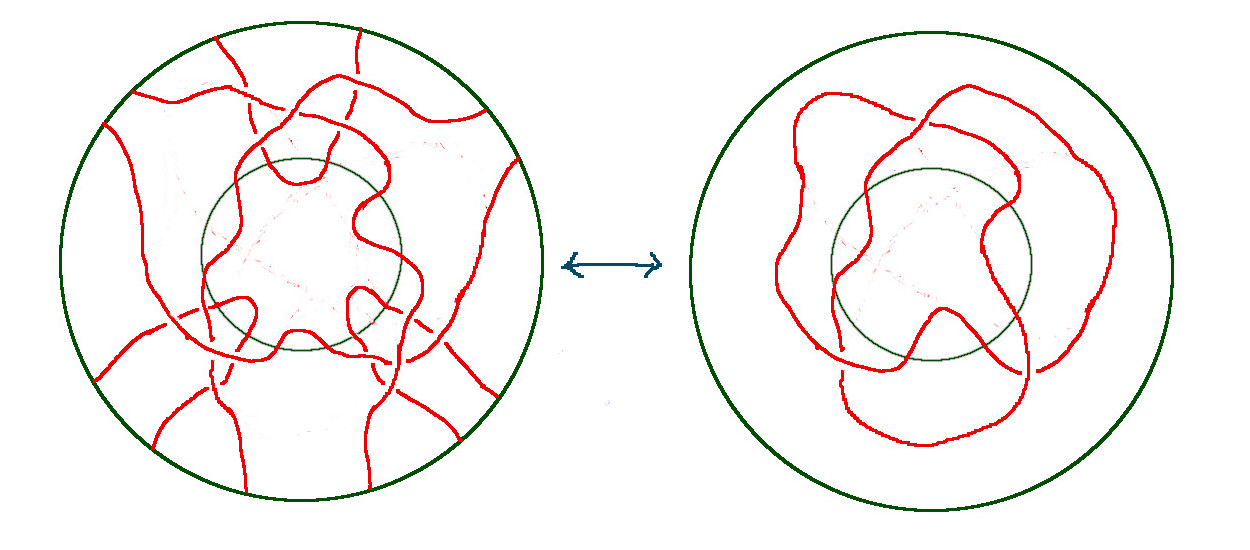}
\caption{A plat representation for affine trefoil}
\end{figure}
\begin{figure}
\includegraphics{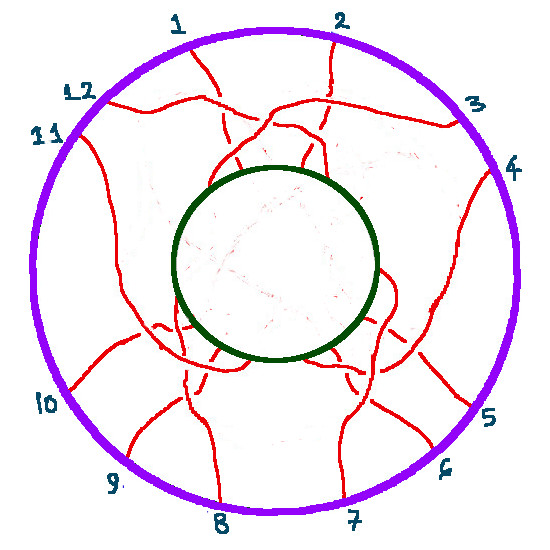}
\caption{The corresponding braid}
\end{figure}

\section{The braid group of $S^2$}
It is easy to see that the composition of motions as defined in section 2, defines a product of braids of the 2-sphere. There is always a \say{identity braid}, which corresponds to all points being at rest. For every motion, there is a \say{inverse} motion which is obtained by reversing the direction of time. Now if we consider the set of all motions of a fixed finite set, say of $n$ points, there is an equivalence relation induced by ambient isotopy relative to the boundary of the strip. The set of isotopy classes clearly forms a group, which is called as the braid group of the sphere. The classical Artin braid group may be described similarly as the braid group of the plane. Refer to \cite{birman}. 

The braid group of $S^2$ is very similar \cite{birman}, to the Artin braid group, which is the braid group of $\mathbb{R}^2$. Here for the sake of simplicity of notation, we will denote the $n-$string braid group of $S^2$ as $B_n$. It should not be confused with the Artin braid group. Since for the purposes of this paper, the only braids we use are from the braid group of $S^2$, there is no chance of confusion. \\
\par Let $C$ be an equator for $S^2$ as before.  We may assume that the boundary points of each string lie on the boundary circles of $C\times I$. That is we are thinking of every braid as motion of finitely many special points on $C$. If $p_1,p_2,...,p_n$ are these special points on $C$, we number both the points $(p_i,0)$ and $(p_i,1)$ by $i$. We choose to index them in a clock wise order. Also it is helpful to think of the indices as elements of $\frac{\mathbb{Z}}{n\mathbb{Z}}$. For keeping the notation simple, we will denote the class of a number say $i+n\mathbb{Z}$, also as just $i$. See Figure 9. \\
\par We can describe the generators of $B_n$ as follows. Consider the braid formed by a crossing between the $i^{th}$ string and the $i+1^{th}$ string and connecting every other special point to the other point with the same index in the obvious way. Refer to Figure 12. Clearly there are two such braids as shown in the diagram and they are inverses of each other in $B_n$. We denote them as $\sigma_i$ and $\sigma_i^{-1}$. Notice that since we have $n+1=1 mod(n)$, we also have $\sigma_n$ and $\sigma_n^{-1}$ which have a crossing between the $n^{th}$ string and the $1^{st}$ string. Clearly we have the following presentation of $B_n$.
\begin{equation*}
B_n\approx  \left\langle \sigma_1,\sigma_2,...,\sigma_n\vert\begin{array}{c}\sigma_i\sigma_j=\sigma_j\sigma_i,\ i-j>2,\\\sigma_i\sigma_{i+1}\sigma_i=\sigma_{i+1}\sigma_i\sigma_{i+1},\\\sigma_1\sigma_2...\sigma_{n-2}\sigma_{n-1}^2\sigma_{n-2}...\sigma_2\sigma_1=1,\\\sigma_{n}=\sigma_1^{-1}\sigma_2^{-1}...\sigma_{n-2}^{-1}=\sigma_{n-1}^{-1}\sigma_{n-2}^{-1}...\sigma_2^{-1}.\end{array}\right\rangle
\end{equation*}
  
\begin{figure}
\includegraphics{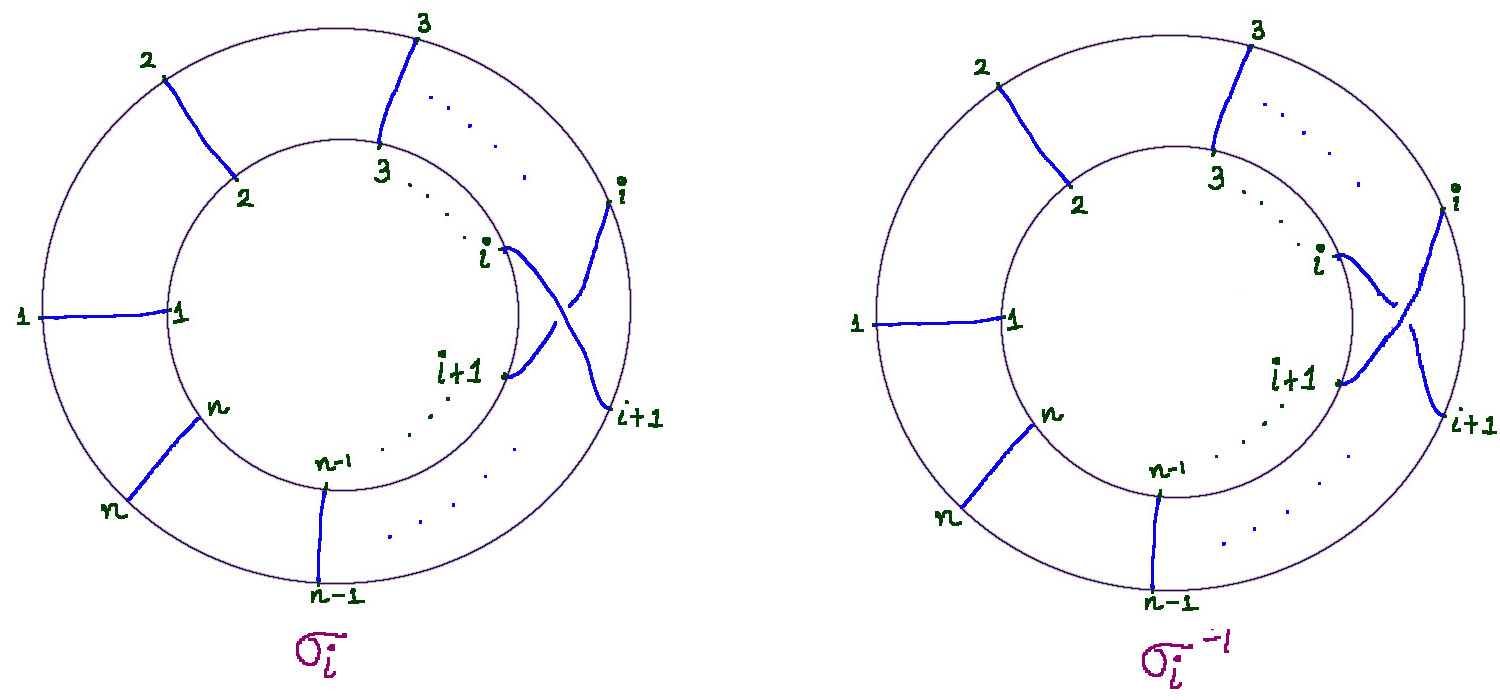}
\caption{Generating braids of $B_n$}
\end{figure}

\section{Residual permutations}

\par There are certain natural questions one would like to ask about the plat representations of links in projective space. For example, by looking at the braid, can we predict the number of components of the closure link? We try to answer this question here. \\
\par The indexing on the boundary points of a braid, $\beta\in B_{k}$, gives a bijection,
\begin{equation*}
f_\beta: \faktor{\mathbb{Z}}{k\mathbb{Z}}\to \faktor{\mathbb{Z}}{k\mathbb{Z}},
\end{equation*}
which we choose to be the one sending the indices of points on $C\times\{0\}$ to the indices of points on $C\times\{1\}$. This is also the projective analogue of the permutation assosiated to a classical braid. We also consider another permutaion,
\begin{equation*}
g: \faktor{\mathbb{Z}}{k\mathbb{Z}}\to \faktor{\mathbb{Z}}{k\mathbb{Z}}
\end{equation*}
defined as follows,
\begin{equation*}
\begin{aligned}
g(i)&= i+1 \text{, if i is an odd class},\\
    &= i-1 \text{, if i is an even class}.
\end{aligned}
\end{equation*}
Since $k=2n$ is even, this is a well defined permutation of $\faktor{\mathbb{Z}}{k\mathbb{Z}}$. Let $G$ denote the group $\faktor{\mathbb{Z}}{k\mathbb{Z}}$.  Notice that $n$ is an order $2$ element in $G$ and let $H$ denote the subgroup generated by $n$. Then we have, 
\begin{equation*}
\faktor{G}{H}\approx \faktor{\mathbb{Z}}{n\mathbb{Z}}.
\end{equation*}
For brevity of notation, we will denote the point $p_i$ on both $C\times\{0\}$ and $C\times\{1\}$ by simply $i$. We may assume that, the points on both $C\times\{0\}$ and $C\times\{1\}$ are arranged symmetrically like numbers on a clock. Then the points $i$ and $i+n$ are diametrically opposite. Thus they belong to the same coset in $\faktor{G}{H}$. We will denote by $[i]$ the coset where the element $i$ belongs. Consider the permutation $f_\beta^{-1}gf_\beta$. Notice that this induces a permutation,
\begin{equation*}
h_\beta: \frac{G}{H}\to \frac{G}{H}.
\end{equation*}
We call this the \textbf{residual permutation} of $\beta$. \\
\begin{theorem}
The number of components in the plat closure link of a braid $\beta$ is same as the number of disjoint cycles in its residual permutation. 
\end{theorem}

 \textbf{Proof:} Notice that, the point $f_\beta^{-1}gf_\beta(i)$ is connected to the point $i$ by the arc formed by the string of $\beta$ connecting $i$ to $f_\beta(i)$ followed by the string in the internal tangle connecting $f_\beta(i)$ to $g(f_\beta(i))$ and then by the string of beta connecting $g(f_\beta(i))$ to $f_\beta^{-1}(g(f_\beta(i)))$. Also notice that each coset in $\faktor{G}{H}$ has two points on $C\times\{0\}$ which are connected by one string in the residual tangle.\\
\par Now suppose $([i_1]\ [i_2]\ ...\ [i_l])$ is a disjoint cycle decomposition of $h_\beta$. Choose any element $j_1$ in the coset $[i_1]$. Then $j_1$ and $j_2:= f_\beta^{-1}(g(f_\beta(j_1)))$ by an arc as described above. And $j_2$ is an element of $[i_2]$ by definition. We can follow this arc from $j_2$ through the string of residual tangle to the point $j_2':= j_2+n \in [i_2]$. If $[i_2]=[i_1]$, in which case $h_\beta$ fixes the $[i_1]$, then this $j_2'=j_1$ and we have a closed loop. Otherwise, we may start again from the point by the above method to the point $j_3:= f_\beta^{-1}(g(f_\beta(j_2')))$. And again by following the string of the residual tangle, we can reach the point $j_3':=j_3+n$. Then, if $j_3'= j_1$, we have a closed loop and the cycle was a transposition $([i_1]\ [i_2])$. Similarly by following this procedure we can see that when we start from the point $j_1$ and reach the end of the cycle at the point $j_l'$, we obtain a knot in the closure link. It is easy to see that if we had chosen to begin at the other element $j_1+n\in [i_1]$ then we would be moving through the same knot, but in the opposite direction. Also if we had represented the cycle with a different ordering of points, then we will again be on the same knot, but starting and ending at a different point.\\
\par Thus every cycle in the disjoint cycle decomposition of $h_\beta$ corresponds to a knot in the plat closure of $\beta$. That is, disjoint cycles in the residual permutation of $\beta$ and the knots in the plat closure of $\beta$ are in one to one correspondence. Hence we are done. \qed

\begin{figure}
\includegraphics{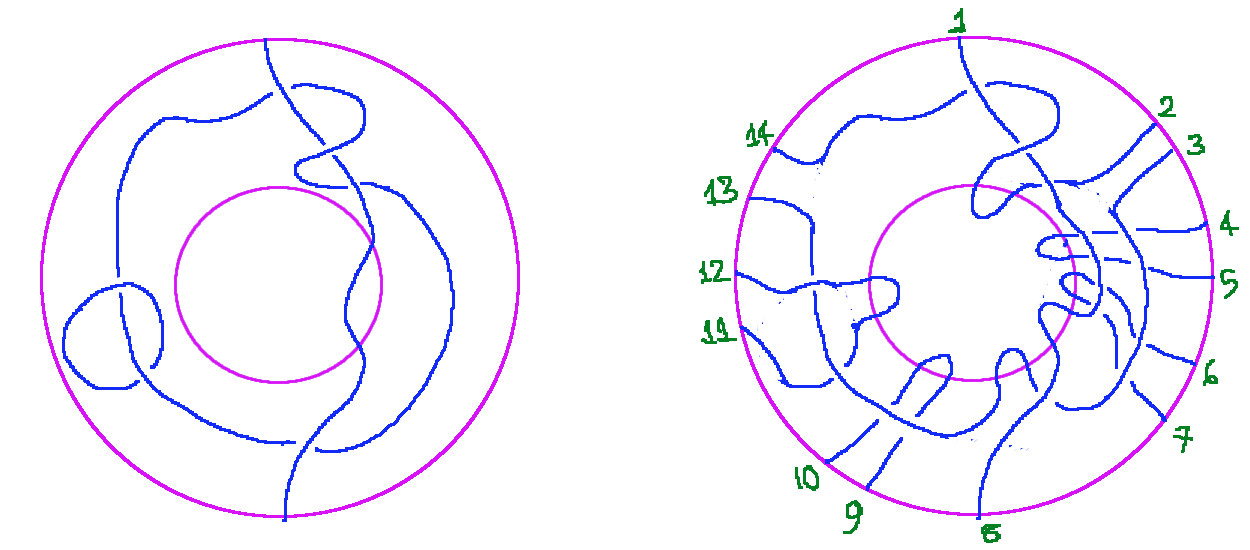}
\caption{The corresponding permutation: (1)(2 7 6 3)(4 5)}
\end{figure}

\par Another natural question, one would like to ask is about the nature of the link formed by closing a braid. Like, its homological properties, affineness and so on. The following theorem studies the conditions for the link to be affine. 
\begin{theorem}
Let $L$ be the closure of an $k=2n$ braid $\beta=\sigma_{i_1}\sigma_{i_2}...\sigma_{i_l}$. Then $L$ is affine if and only if there exists an even class $j \in \faktor{\mathbb{Z}}{k\mathbb{Z}}$ such that, both $\sigma_j$ and $\sigma_{j+n}$ don't appear in the word representing $\beta$. 
\end{theorem}
\textbf{Proof:}Suppose $j$ is such an even class. Then in a diagram of the closure link, drawn on a disk, we may draw a diametrical line on the diagram which passes through the region of the strip between the points $j$ and $j+1$. The line  doesn't intersect the link. And hence, by pulling it back under the projection, we can construct a projective plane disjoint from the link. Thus the link is affine, by Theorem 1.2 in \cite{ramavis}.\\ 
\par Now we will prove the converse. If $L$ is affine, choose a projection, with a dijoint diameter $d$ for the disk. There are four points in $d\cap C\times\{0,1\}$, two on each of the circles. Then, in the strip, we can choose an indexing for the boundary points of $\beta$ so that, when we move clock wise from the an intersection point on $C\times\{0\}$, the first on the boundary point of $\beta$ we reach is indexed as $1$. Notice that the choice of this intersection point is arbitrary, but any of the two points in $d\cap C\times\{0\}$ will work. Then, the boundary point we see first in the anti-clockwise direction on $C\times\{0\}$ is indexed as $2n$, an even class. This indexing induces an indexing on $C\times\{1\}$. Then choose, $j=2n$, then it will satisfy the condition above. Hence we are done. \qed\\

\section{Moves on braids}
As usual in the classical case, the same link may be represented by plat closures of  multiple braids. We will now explore the equivalence relation between the braids which have isotopic closures.\\

\begin{figure}
\includegraphics{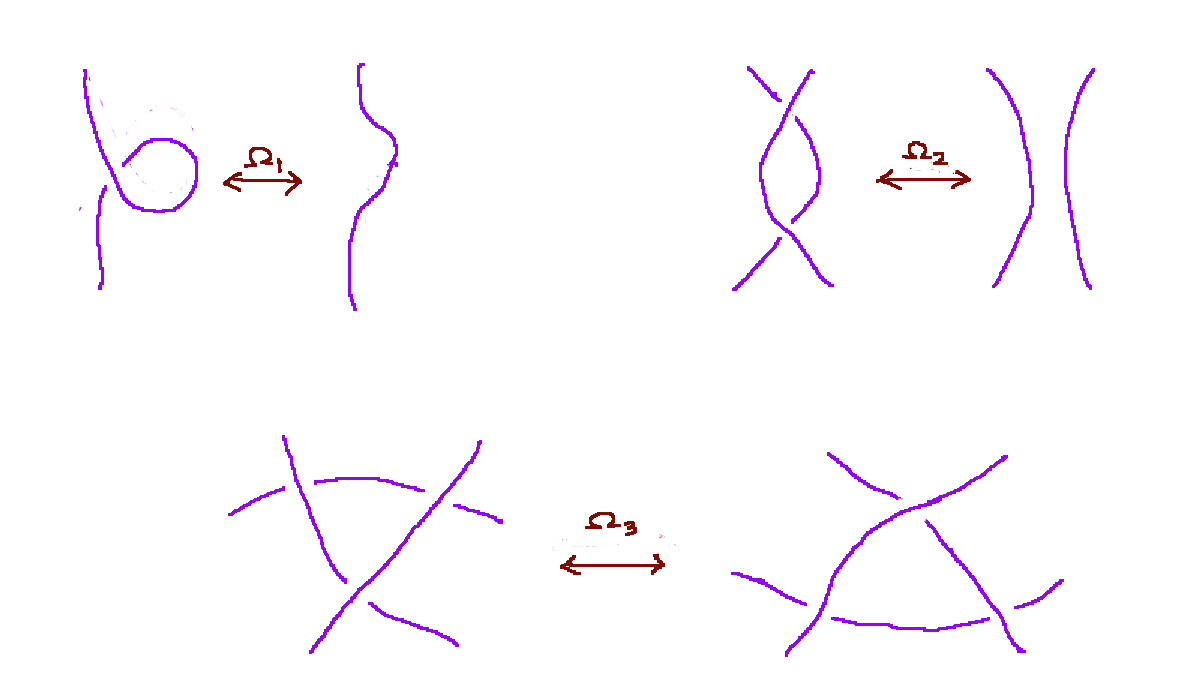}\\
\includegraphics{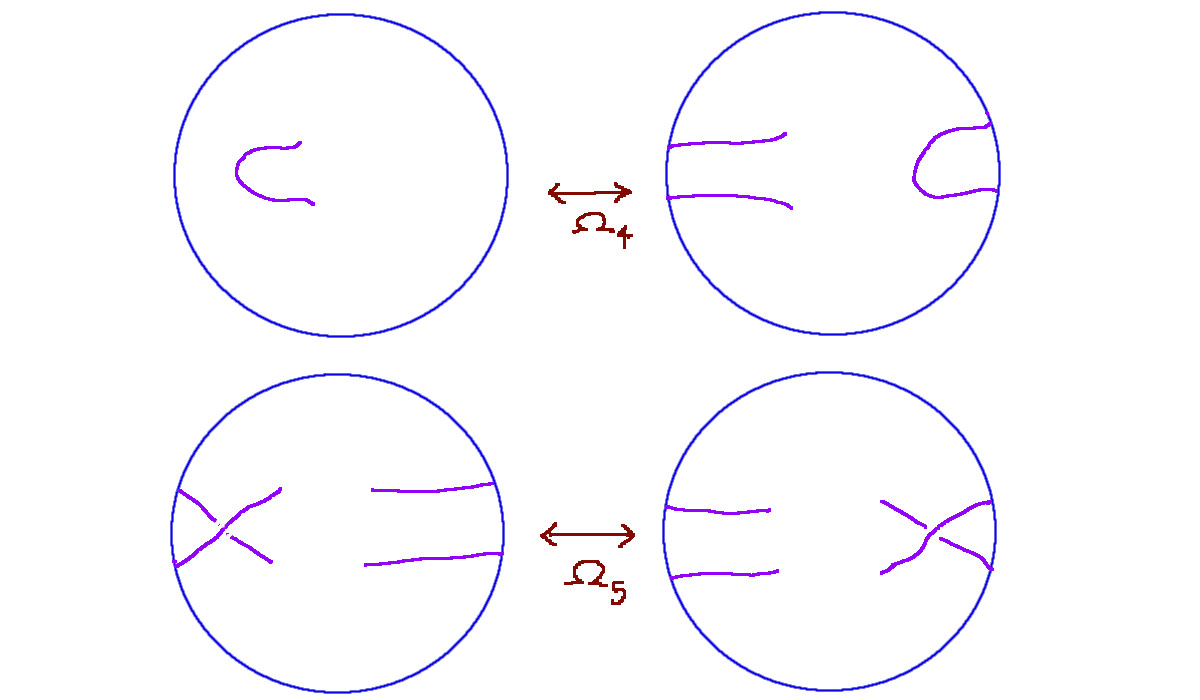}
\caption{Moves of projective diagrams}
\end{figure}

The diagramatic moves in $\mathbb{R}P^3$, as described in \cite{julia} are as shown in Figure 14. \textit{Two links in $\mathbb{R}P^3$ are isotopic if and only if their diagrams can be transformed from one to the other by a finite sequence of the moves $\Omega_1, \Omega_2, \Omega_3, \Omega_4$ and $\Omega_5$.}
That is, as sets, the set of equivalence classes of diagrams induced by above described moves is in one to one correspondence with the set of isotopy classes of links.\\
Suppose $\beta:=\sigma_{i_1}\sigma_{i_2}...\sigma_{i_k}$ is a braid. Notice that the move $\Omega_2$ applied on $\beta$ is equivalent to deleting a pair of the form $\sigma_{i_l}\sigma_{i_l}^{-1}$ from the word representing $\beta$. And performing a $\Omega_3$ move on $\beta$ is equivalent replacing a sub-word of the form $\sigma_{i}\sigma_{i+1}\sigma_{i}$ with a word $\sigma_{i+1}\sigma_{i}\sigma_{i+1}$ for an arbitrary $i$. Now we know this is already a relation in the braid group $B_k$. Clearly the inverse of both these processes can also be described in a similar manner. That is, the moves $\Omega_2$ and $\Omega_3$ are already incorporated in the group structure of the braid group.  \\
\par Consider the following moves on braids. Suppose $k=2n$ is even. Let $\beta=\sigma_{i_1}\sigma_{i_2}...\sigma_{i_m}$ be a braid in $B_k$. Then there will be exactly $n$ odd classes, i.e, $(2l+1)+k\mathbb{Z}$ in $\faktor{\mathbb{Z}}{k\mathbb{Z}}$. \\
\par Let $M_0^i$ be a move as follows,
\begin{align*}
M_0^i: \beta\sigma_{i-1}\sigma_i \longleftrightarrow \beta\sigma_{i-1}^{-1}\sigma_i^{-1}, \text{if $i$ is odd}.\\
\overline{M_0^i}: \beta\sigma_{i-1}^{-1}\sigma_i^{-1} \longleftrightarrow \beta\sigma_{i-1}\sigma_i, \text{if $i$ is even}.
\end{align*}
\begin{figure}
\includegraphics[scale=0.75]{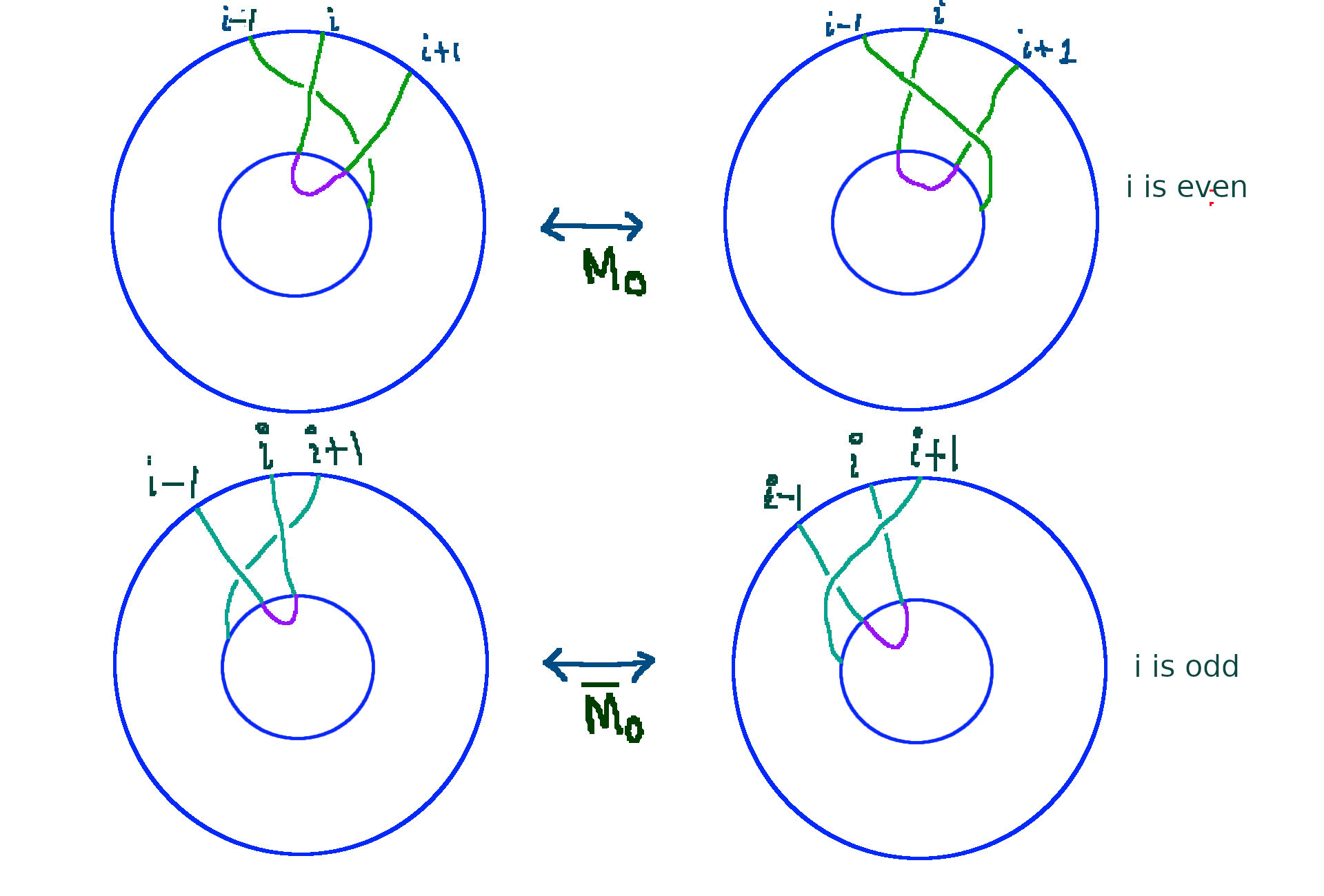}
\caption{$M_0^i$ moves}
\end{figure}
Figure 15 demonstrates this move. It is easy to see that the plat closures of two braids which are related by a move of this type are equivalent by two $\Omega_2$-moves. 
\par For an odd $i$, let $M_1^i$ be the move,
\begin{align*}
M_1^i: \beta\longleftrightarrow \beta\sigma_i,\\
\overline{M_1^i}: \beta\longleftrightarrow \beta\sigma_i^{-1}. 
\end{align*}

\begin{figure}
\includegraphics[scale=0.55]{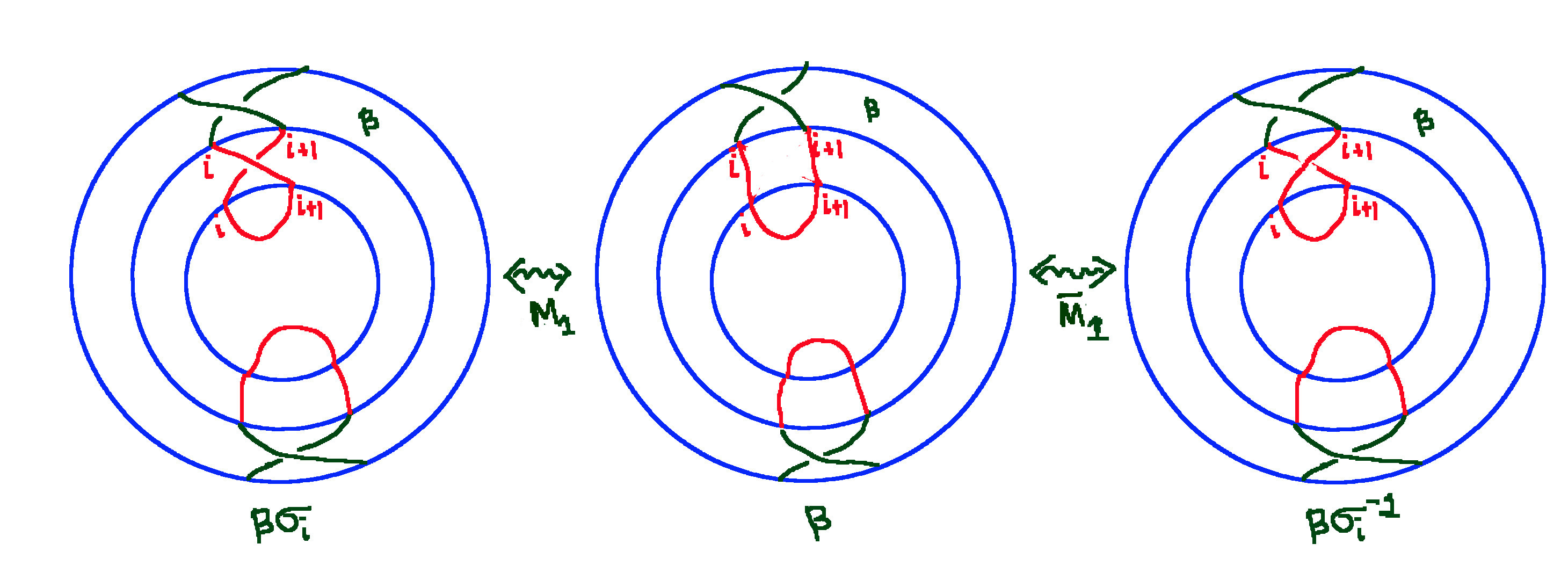}
\caption{A typical $M_1$-move}
\end{figure}

Notice that the internal tangle for a plat has strings connecting every odd class $i$ with $i+1$. Then the diagrams of plat closures of two braids  which are related by an $M_1^i$ move are related by an $\Omega_1$ move. Thus the braids have isotopic plat closures. Refer to Figure 16. \\
\par Now suppose $\beta=\sigma_l\beta'$ is a braid in $B_k$ where $k=2n$ as above. Then we have the following move on the $B_k$.
\begin{align*}
M_2: \sigma_l\beta'\longleftrightarrow \beta'\sigma_{l+n}\\
\overline{M_2}: \sigma_l^{-1}\beta'\longleftrightarrow \beta'\sigma_{l+n}^{-1}
\end{align*} 

\begin{figure}
\includegraphics[scale=0.45]{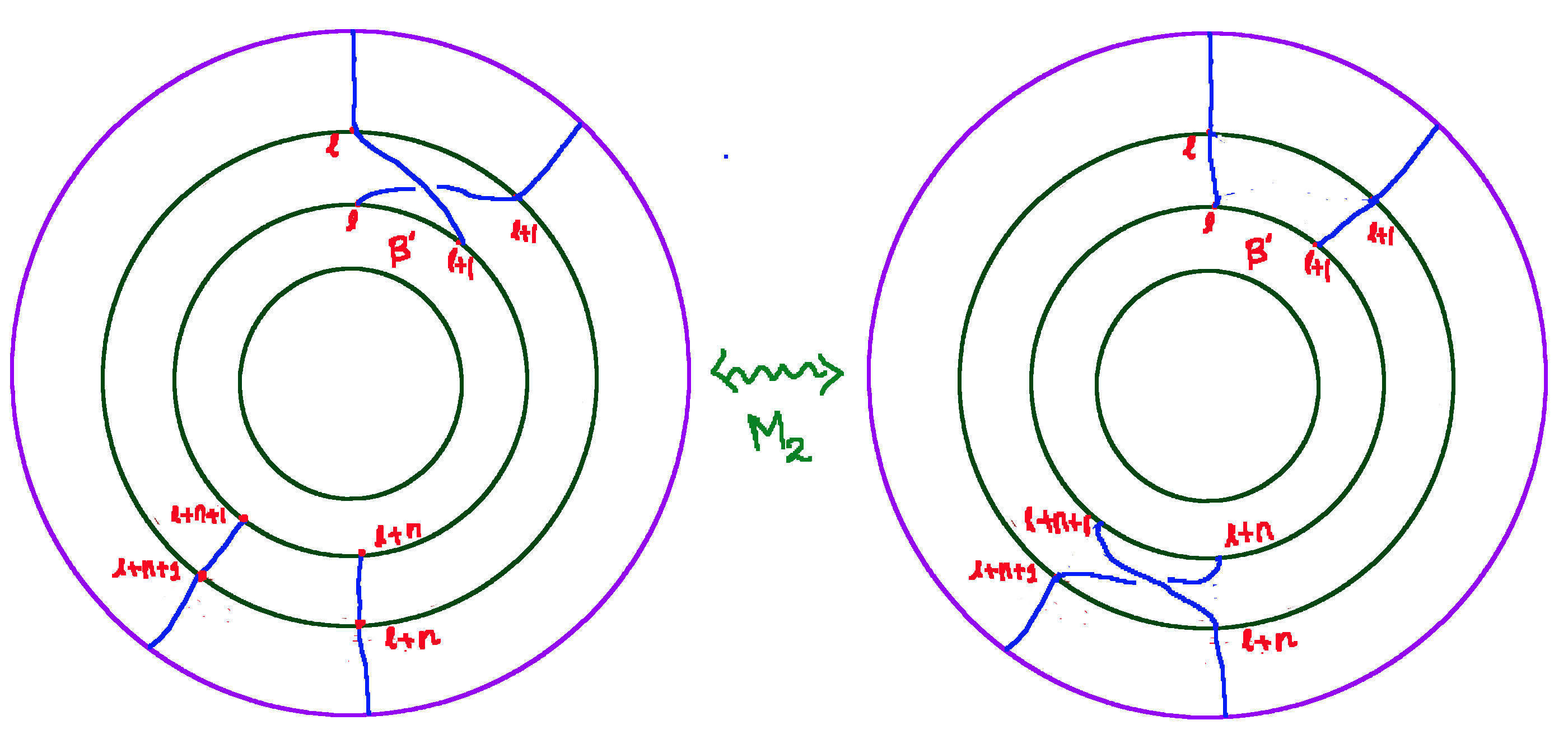}\\
\includegraphics[scale=0.45]{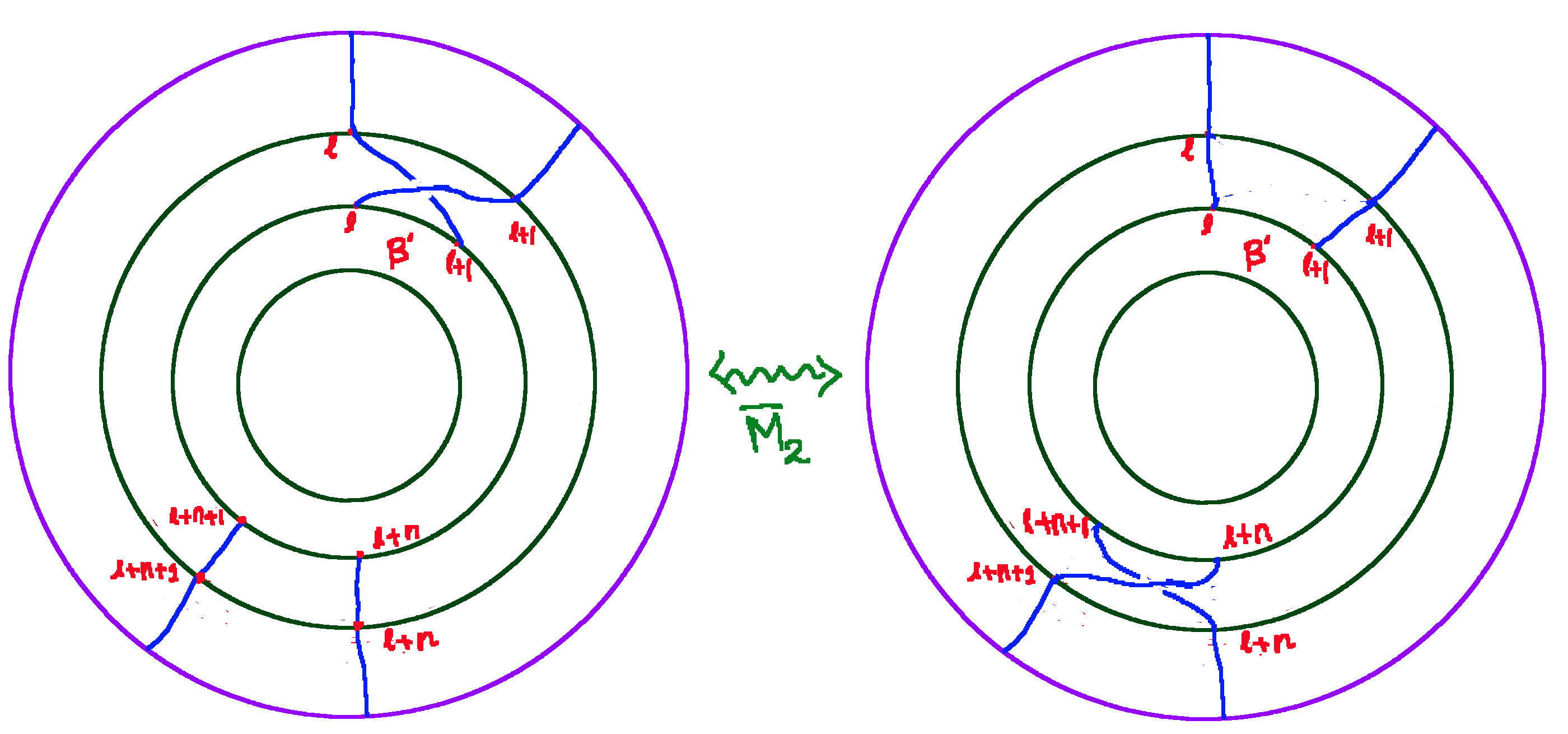}
\caption{$M_2$-moves}
\end{figure}
Clearly the plat closures of two braids related by an of $M_2$ or an $\overline{M_2}$ moves, are isotopic since their diagrams are equivalent through an $\Omega_5$ move. \\
\par Suppose we think of the braid group $B_{k+4}$ as the group of motion of the points $p_1,p_2,...,p_n, p_{n+1}, p_{n+2},-p_1,-p_2,...,-p_n, -p_{n+1}, -p_{n+2}$ on $C$ where each $p_i$ is next to $p_{i+1}$ for $i<n+2$. Let $\gamma$ be a braid formed by a motion where the points $p_1,p_2,-p_1,-p_2$ were still. Then each of the strings formed by these four points in $\gamma$, is not braided with any other string in $\gamma$. If we remove these four strings from $\gamma$ it would not disturb the motion of the other $k$ points which will form a $k$ braid. Thus we may think of each of the $k+4$ string braids formed by a motion where the points $p_1,p_2,-p_1,-p_2$ are still as obtained from a $k$ string braid by introducing four still points on $C$. Thus we have natural map,
\begin{equation*}
e: B_k\hookrightarrow B_{k+4}
\end{equation*}

\begin{figure}
\includegraphics{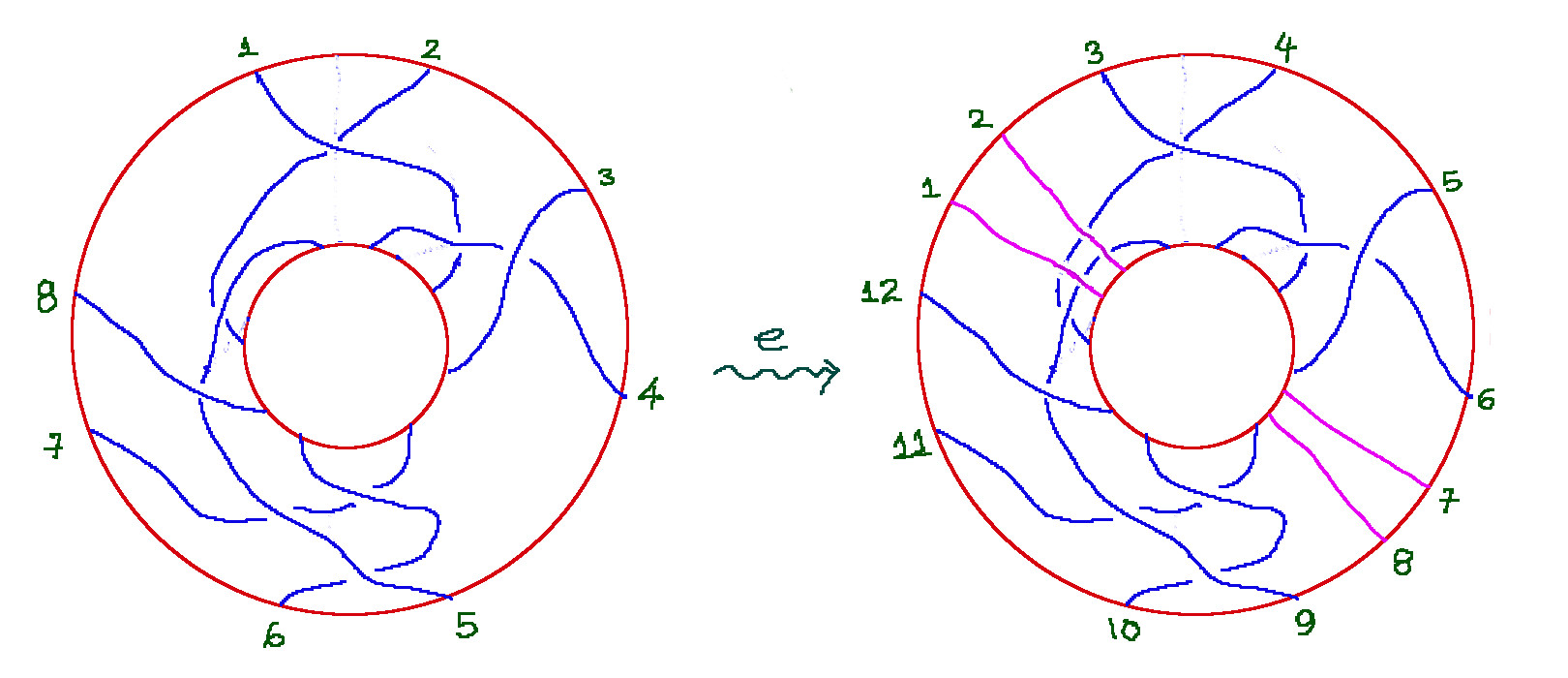}
\caption{A typical case of applying the map $e: B_8\hookrightarrow B_{12}$}
\end{figure}

induced by introducing four points on the equator $C$ of the type $x$, $y$, $-x$, $-y$ such that no points of $B_k$ are lying between $x$ and $y$. See Figure 18. Suppose the braids in $B_k$ were described as the motion of $q_1,q_2,...,q_n, -q_1,-q_2,...,-q_n$. The image of each $k$-braid under $e$  is a $k+4$ braid with the four more strings $\{p_1\}\times I$, $\{p_2\}\times I$, $\{-p_1\}\times I$ and $\{p_1\}\times I$ in the strip. We choose to rename the points in the strip in such a way that, $x$, $y$, $-x$ and $-y$ will be relabeled as $p_1$, $p_2$, $-p_1$ and $-p_2$. Then the rest of the points will be relabeled according to the ordering. That is $q_1,q_2,...,q_n$ will be labeled as $p_3,p_4,...p_{n+2}$ and $-q_1,-q_2,...,-q_n$ as $-p_3, -p_4,..., -p_{n+2}$ respectively.  \\
\par For any $3\leq l\leq k+2$, define 
\begin{align*}
\alpha_l:=\sigma_2\sigma_3\sigma_4...\sigma_{l-2}\sigma_{l-1}\sigma_{l-2}^{-1}\sigma_{l-3}^{-1}...\sigma_2^{-1},\\
\overline{\alpha_l}:=\sigma_2^{-1}\sigma_3^{-1}\sigma_4^{-1}...\sigma_{l-2}^{-1}\sigma_{l-1}^{-1}\sigma_{l-2}\sigma_{l-3}...\sigma_2.
\end{align*}
certain family braids in $B_{k+4}$. 

\begin{figure}
\includegraphics[scale=0.75]{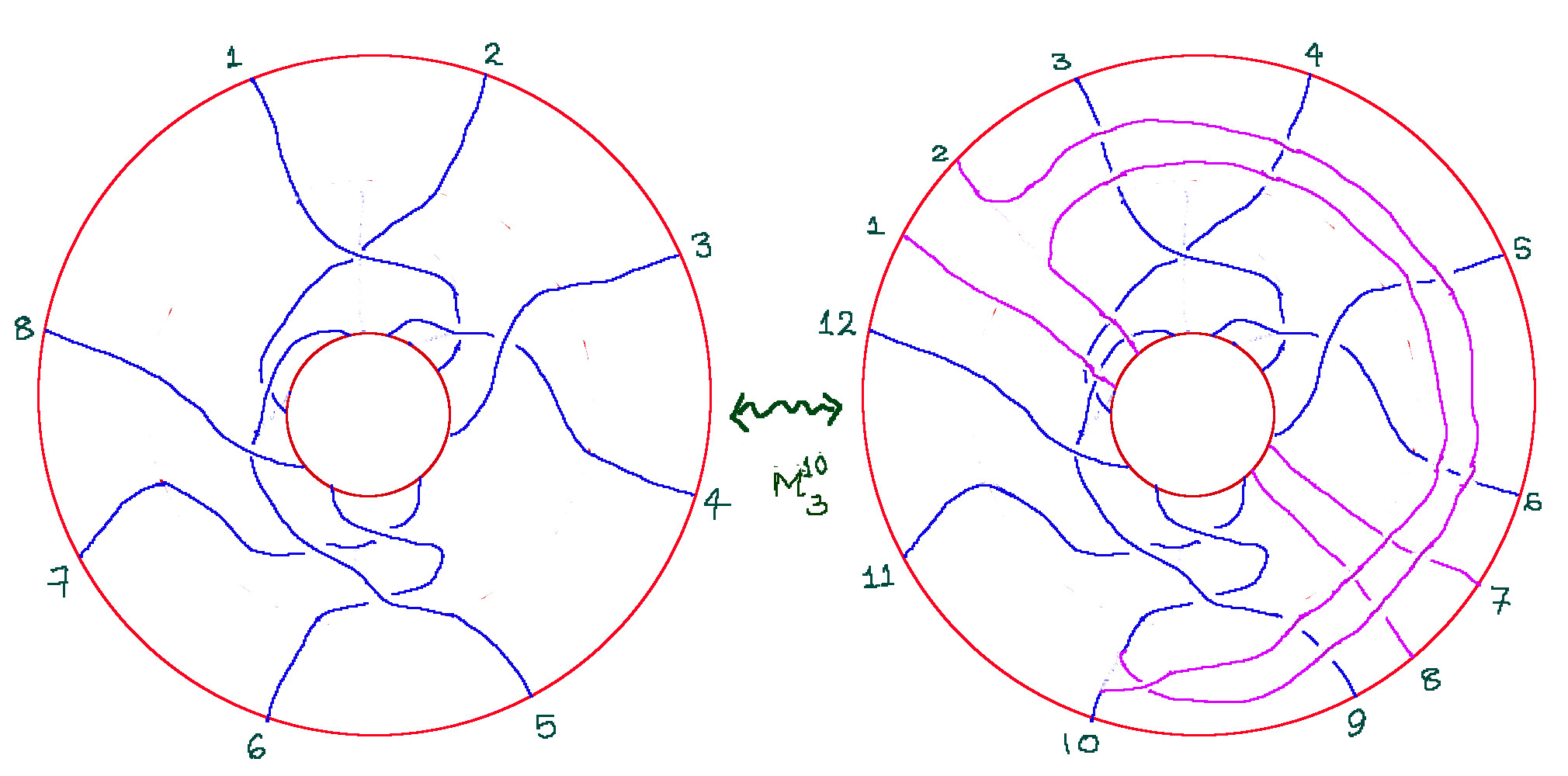}
\caption{$M_3^{10}$-move performed on a braid in $B_8$ resulting in $B_{12}$}
\end{figure}

 Then we define,
\begin{align*}
M_3^{l}: \beta \longleftrightarrow \alpha_l e(\beta),\\
\overline{M_3}^l: \beta \longleftrightarrow \overline{
\alpha_l} e(\beta).
\end{align*}
Which is a move relating braids of $B_k$ and $B_{k+4}$. Refer to Figure 19. Each instance of performing the above move and obtaining a braid $\beta'$ from a $\beta\in B_k$, the diagram of their plat closure of $\beta$ change by an $\Omega_4$ move after an $\Omega_1$ move. Thus clearly, the plat closures of braids which are related by these moves are isotopic.\\

There is another move which also changes the index of the braid. Suppose one string in a braid is isotopically moved to form a pair of maxima and minima of the projection, $f:S^2\times I\to I$. Notice that then its no longer monotonic on this string. We may bring the extrema points into the ball region by isotopically moving them. Figure 20 and 21 presents the cases when half of the braid index is odd and even. We call this as $M_4$ and $\overline{M_4}$ moves. When $n$ is odd, the move is just an application of the map, $e:B_k\to B_{k+4}$ defined above. For an even $n$, $\overline{M_4}$ clealy defines a map from $B_k$ to $B_{k+4}$. But if $\beta$ is a $k$-braid, the form of the braid $M_4(\beta)$ depends on the form of the braid $\beta$ and writing a closed expression may not be possible, just like the map $e$ mentioned above.  

\begin{figure}
\includegraphics{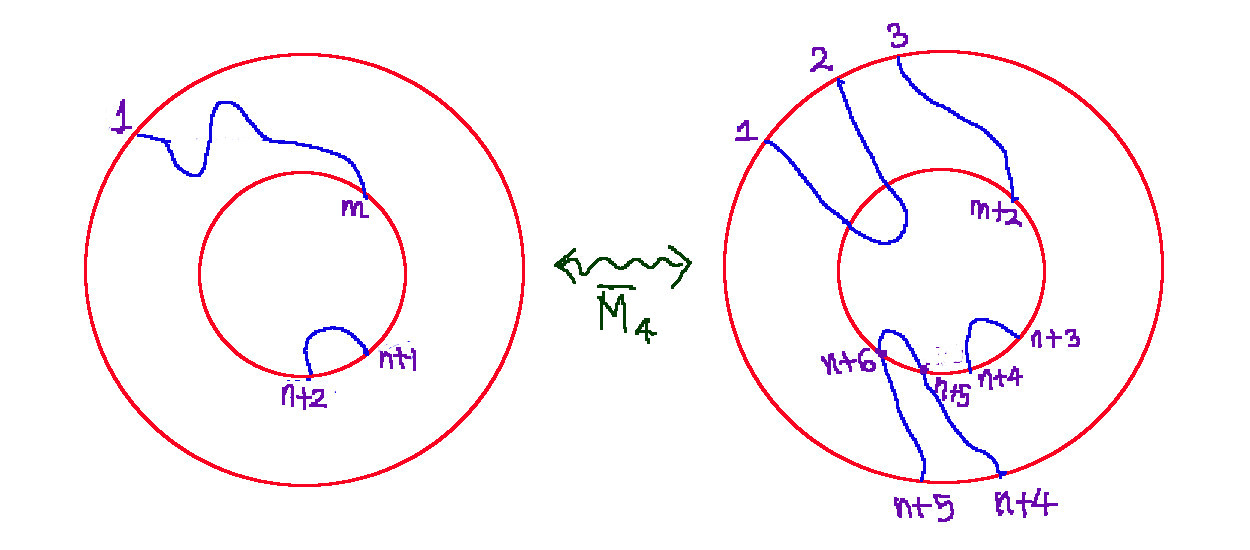}
\caption{When $n$ is even}
\end{figure}
\begin{figure}
\includegraphics{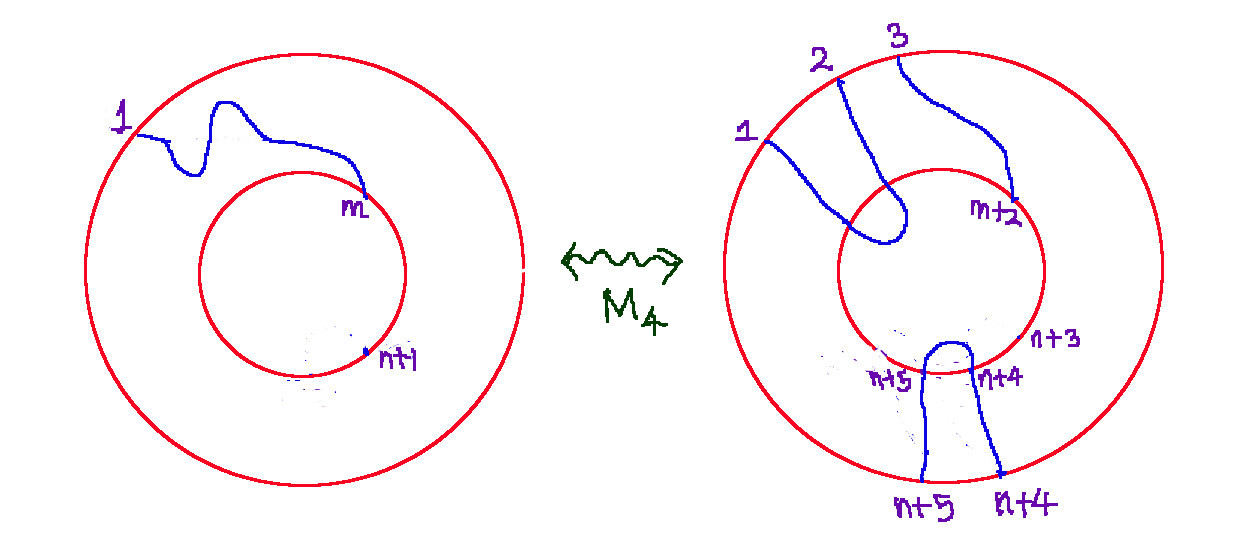}
\caption{When $n$ is odd}
\end{figure}

\par Notice that all the moves described above with a name of the type {\bf $M_i$} always comes with a pair $\overline{M_i}$. In each case, it is obvious to see that each one is just another version of the other. In what follows we will drop the overlines and call these moves as just $M_i$ since all what we are saying is applicable to both with some obvious modifications. We will call the set of all the operations as just, $M$-moves. If a braid $\beta$ can be turned into another braid $\beta'$ by a finite sequence of $M$-moves, we will say $\beta$ and $\beta'$ are \textbf{$M$-equivalent}. In the above discussions we observe that plat closures of $M$-equivalent braids are isotopic links in $\mathbb{R}P^3$. We also propose the following conjecture.\\

\begin{conjecture}
The projective plat closures of two braids are isotopic if and only if they are $M$-equivalent. 
\end{conjecture}

\par As a support for the conjecture we provide Figure 22 where we show that the braids in two different plat closure diagrams of the figure-8 knot are $M$-equivalent.\\

\begin{figure}
\includegraphics{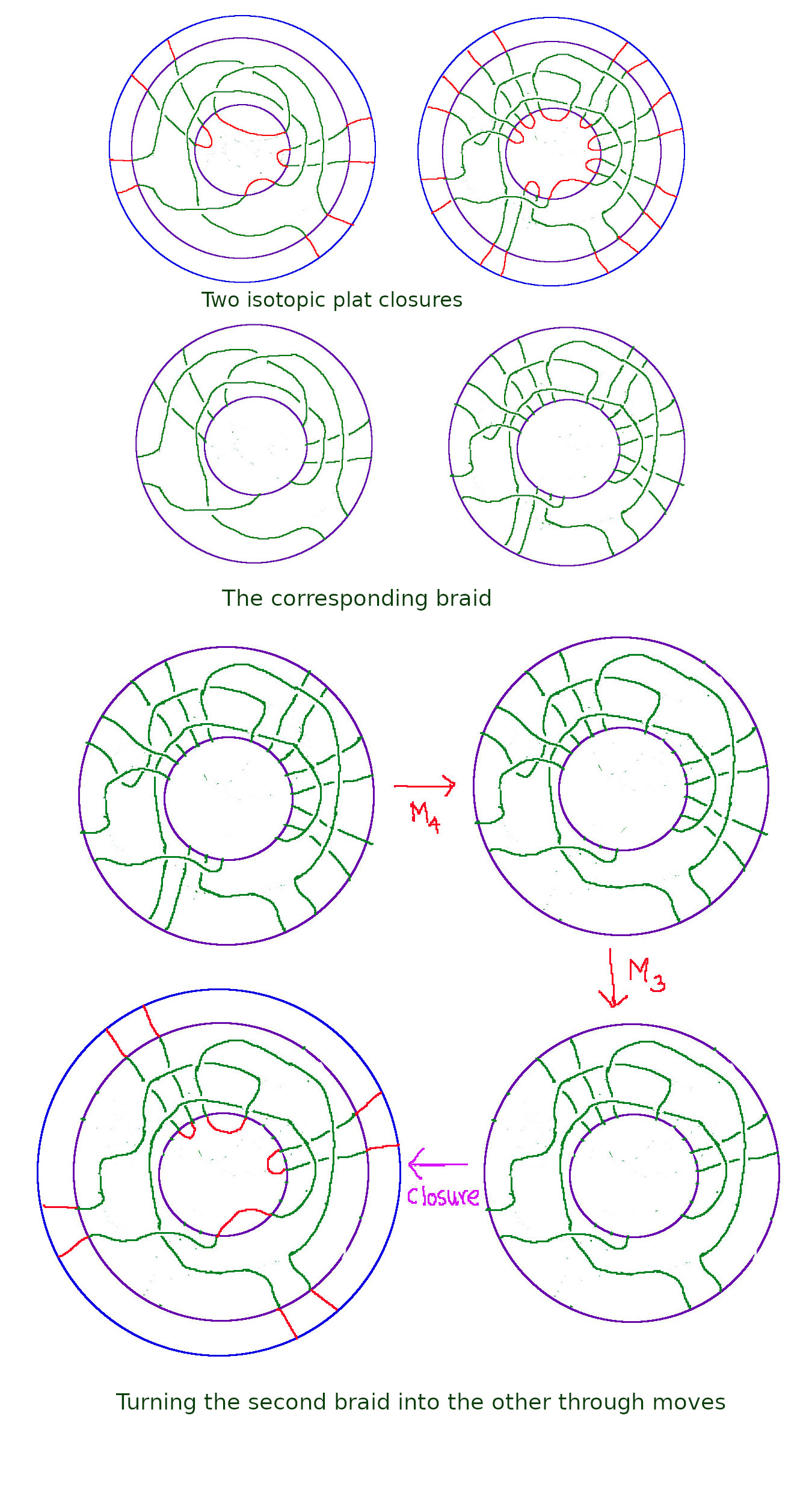}
\caption{An example of applying the moves}
\end{figure}

\textbf{Concluding remark:} One of the wonderful applications of this work, may possibly be that one can use it for constructing some quantum invariants for knots and links in $\mathbb{R}P^3$. This can be done by finding an appropriate series of representations, $V_k$ for $B_k$ which admits a good \say{quantum trace}, $t_k: B_k \to R$, (where $R$ is some appropriate ring, for example $\mathbb{C}[x,x^{-1}]$), invariant under $M$-moves. This trace will be a link invariant.


\begin{thebibliography}{9}

\bibitem{birman}
Birman, J. S., \textit{Braids, Links, and Mapping Class Groups}.(AM-82), Volume 82. Vol. 82. Princeton University Press, 2016.

\bibitem{alex}Alexander, J. W., \textit{A lemma on systems of knotted curves.}  Proceedings of the National Academy of Sciences, 9(3), 93-95 (1923)


\bibitem{julia}
Drobotukhina, Y. V., \textit{An analogue of the Jones polynomial for links in RP3 and a generalization of the Kauffman-Murasugi theorem.}, Leningrad Math. J. 2:3 : 613-630 ,1991

\bibitem{ramavis} Mishra R., Narayanan V., \textit{Geometry of knots in real projective $3$-space}, (To appear in) \textit{Journal of Knot Theory and its Ramifications.}

\end{thebibliography}
\end{document}